# LAW OF THE ITERATED LOGARITHM
# FOR STATIONARY PROCESSES[1]

By Ou Zhao and Michael Woodroofe

*University of Michigan*

There has been recent interest in the conditional central limit question for (strictly) stationary, ergodic processes $\ldots, X_{-1}, X_0, X_1, \ldots$ whose partial sums $S_n = X_1 + \cdots + X_n$ are of the form $S_n = M_n + R_n$, where $M_n$ is a square integrable martingale with stationary increments and $R_n$ is a remainder term for which $E(R_n^2) = o(n)$. Here we explore the law of the iterated logarithm (LIL) for the same class of processes. Letting $\|\cdot\|$ denote the norm in $L^2(P)$, a sufficient condition for the partial process to have the form $S_n = M_n + R_n$ is that $n^{-3/2}\|E(S_n|X_0, X_{-1}, \ldots)\|$ be summable. A sufficient condition for the LIL is only slightly stronger, requiring $n^{-3/2} \log^{3/2}(n)\|E(S_n|X_0, X_{-1}, \ldots)\|$ to be summable. As a by-product of our main result, we obtain an improved statement of the conditional central limit theorem. Invariance principles are obtained as well.

**1. Introduction.** Let $\ldots, X_{-1}, X_0, X_1, \ldots$ denote a centered, square integrable, (strictly) stationary and ergodic process, defined on a probability space $(\Omega, \mathcal{A}, P)$, with partial sums denoted by $S_n = X_1 + \cdots + X_n$. The main question addressed is the law of the iterated logarithm: under what conditions is

$$(1) \qquad \limsup_{n \to \infty} \frac{S_n}{\sqrt{2n \log_2(n)}} = \sigma \qquad \text{w.p. } 1$$

for some $0 \le \sigma < \infty$, where $\log_2(n) = \log(\log(n))$. Of course, (1) holds if the $X_i$ are independent, by the classic work of Hartman and Wintner [6], and more generally—for example, [7, 15, 17]. Here we employ an approach which has been used recently in the study of the central limit question for stationary processes—martingale approximations.

Received March 2006; revised October 2006.

[1]Supported by the National Science Foundation.

*AMS 2000 subject classifications.* Primary 60F15; secondary 60F05.

*Key words and phrases.* Conditional central limit question, ergodic theorem, Fourier series, martingales, Markov chains, operators on $L^2$.

As in Maxwell and Woodroofe [11], it is convenient to suppose that $X_k$ is of the form $X_k = g(W_k)$, where $\ldots, W_{-1}, W_0, W_1, \ldots$ is a stationary, ergodic Markov chain. The state space, transition function and (common) marginal distribution are denoted by $\mathcal{W}, Q$ and $\pi$; thus, $\pi(B) = P[X_n \in B]$, and

$$Qf(w) = E[f(W_{n+1})|W_n = w]$$

for a.e. $w \in \mathcal{W}$, measurable $B \subseteq \mathcal{W}$ and $f \in L^1(\pi)$. The iterates of $Q$ are denoted by $Q^k$. It is also convenient to suppose that the probability space $\Omega$ is endowed with an ergodic, measure-preserving transformation $\theta$ for which $W_k \circ \theta = W_{k+1}$ for all $k$. Neither convenience entails any loss of generality, since we may let the probability space be $\mathbb{R}^{\mathbb{Z}}$, $X_k$ be the coordinate functions, $W_k = (\ldots, X_{k-1}, X_k)$, and $\theta$ be the shift transformation. Some other choices of $W_k$ are considered in the examples.

Let $\|\cdot\|$ denote the norm in $L^2(P)$, $\mathcal{F}_k = \sigma(\ldots, W_{k-1}, W_k)$, and recall the main result of [11]; if

$$(2) \qquad \sum_{n=1}^{\infty} n^{-3/2} \|E(S_n|\mathcal{F}_0)\| < \infty,$$

then

$$(3) \qquad \sigma^2 := \lim_{n\to\infty} \frac{1}{n} E(S_n^2)$$

exists and is finite, and

$$(4) \qquad S_n = M_n + R_n,$$

where $M_n$ is a square integrable martingale with stationary, ergodic increments, and $\|R_n\| = o(\sqrt{n})$. It is shown in [11] that if (2) holds, then the conditional distributions of $S_n/\sqrt{n}$, given $\mathcal{F}_0$, converge *in probability* to the normal distribution with mean 0 and variance $\sigma^2$ (see their Corollary 1). It can also be shown that (2) is *best possible* through Peligrad and Utev [13].

To state the main result of the paper, let $\ell$ be a positive, nondecreasing and slowly varying (at $\infty$) function and let

$$\ell^*(n) = \sum_{j=1}^{n} \frac{1}{j\ell(j)}.$$

THEOREM 1. *If $\ell$ is a positive, slowly varying, nondecreasing function and*

$$(5) \qquad \sum_{n=1}^{\infty} n^{-3/2}\sqrt{\ell(n)}\log(n)\|E(S_n|\mathcal{F}_0)\| < \infty,$$

*then*

$$\lim_{n\to\infty} \frac{R_n}{\sqrt{n\ell^*(n)}} = 0 \qquad w.p.\ 1.$$

COROLLARY 1. *If* (5) *holds with* $\ell(n) = 1 \vee \log(n)$, *then* (1) *holds.*

PROOF. In this case $\ell^*(n) \sim \log_2(n)$, so that $R_n/\sqrt{n \log_2(n)} \to 0$ as $n \to \infty$, and

$$\limsup_{n \to \infty} \frac{S_n}{\sqrt{2n \log_2(n)}} = \limsup_{n \to \infty} \frac{M_n}{\sqrt{2n \log_2(n)}}$$

both w.p. 1. The corollary now follows from the law of the iterated logarithm of martingales; for example, Stout [17]. $\square$

The next corollary strengthens the conclusion of [11] from convergence in probability to convergence w.p. 1, under a slightly stronger hypothesis. Kipnis and Varadhan [8] call this an important question in a closely related context (see their Remark 1.7). Let $F_n$ denote a regular conditional distribution function for $S_n/\sqrt{n}$ given $\mathcal{F}_0$, so that

$$F_n(\omega; z) = P\left[\frac{S_n}{\sqrt{n}} \leq z \Big| \mathcal{F}_0\right](\omega)$$

for $\omega \in \Omega$ and $-\infty < z < \infty$; and let $\Phi_\sigma$ denote the normal distribution with mean 0 and variance $\sigma^2$.

COROLLARY 2. *If* (5) *holds with some* $\ell$ *for which* $1/[n\ell(n)]$ *is summable, then* $F_n(\omega; \cdot)$ *converges weakly to* $\Phi_\sigma$ *for a.e.* $\omega$.

PROOF. Let $G_n$ be a regular conditional distribution for $M_n/\sqrt{n}$ given $\mathcal{F}_0$. Then $G_n(\omega; \cdot)$ converges weakly to $\Phi_\sigma$ for a.e. $\omega$, essentially by the martingale central limit theorem, applied conditionally given $\mathcal{F}_0$. See [11] for the details. Moreover, $P[\lim_{n \to \infty} R_n/\sqrt{n} = 0 | \mathcal{F}_0] = 1$ w.p. 1, since $P[\lim_{n \to \infty} R_n/\sqrt{n} = 0] = 1$, by Theorem 1. The corollary follows easily. $\square$

A major contribution of this paper is to obtain a simple, general sufficient condition (5) for the LIL. Our results differ from those of Arcones [1], for example, by not requiring normality, and those of Rio [15] by not requiring strong mixing. In [10], Lai and Stout have a quite general result for strongly dependent variables. Their results require a condition on the moment-generating function of the delayed partial sums and only cover the upper half of LIL. Yokoyama [18] also uses martingale approximation in a similar setting to ours. His results require a martingale approximation, as in (4), and bounds on higher moments of the remainder term.

The rest of the paper is organized as follows. The proof of Theorem 1 is outlined in Section 2, with supporting details in Sections 3 and 4. Invariance principles are considered in Section 5, and examples in Section 6.

**2. Outline of the proof.** In this section, we give an outline of the proof for the main result. Let

$$(6) \qquad h_\varepsilon = \sum_{k=1}^\infty \frac{Q^{k-1}g}{(1+\varepsilon)^k}$$

and $H_\varepsilon(w_0, w_1) = h_\varepsilon(w_1) - Q h_\varepsilon(w_0)$. Thus $H_\varepsilon \in L^2(\pi_1)$, where $\pi_1$ denotes the joint distribution of $W_0$ and $W_1$. In [11] it is shown that if (2) holds, then $H := \lim_{\varepsilon \downarrow 0} H_\varepsilon$ exists in $L^2(\pi_1)$ and that (4) holds with $M_n = H(W_0, W_1) + \cdots + H(W_{n-1}, W_n)$. Letting $\xi_k = g(W_k) - H(W_{k-1}, W_k)$ leaves

$$(7) \qquad R_n = \sum_{k=1}^n \xi_k = \sum_{k=1}^n \xi_0 \circ \theta^k$$

in (4).

For appropriately chosen $\beta_k \sim c/\sqrt{k^3 \ell(k)}$ [see (12), below], the series

$$(8) \qquad B(z) = \sum_{k=1}^\infty \beta_k z^k$$

converges for all complex $|z| \le 1$, is analytic in $|z| < 1$, $B(1) = 1$, and $|1 - B(z)| > 0$ for $z \ne 1$. Letting $T$ be the operator on $L^2(P)$ defined by $T\eta = \eta \circ \theta$, it is also true that $B(T)$ converges in the operator norm. Thus,

$$(9) \qquad B(T)\eta = \sum_{k=1}^\infty \beta_k T^k \eta = \sum_{k=1}^\infty \beta_k \eta \circ \theta^k.$$

With this notation, there are two main steps to the proof. It is first shown that in (7), $\xi_0 \in [I - B(T)]L^2(P)$, the range of $I - B(T)$, so that $\xi_0 = \eta_0 - B(T)\eta_0$ for some $\eta_0 \in L^2(P)$. It is then shown that for any $\xi \in [I - B(T)]L^2(P)$,

$$\lim_{n \to \infty} \frac{1}{\sqrt{n\ell^*(n)}} \sum_{k=1}^n T^k \xi = 0 \qquad \text{w.p. } 1.$$

The broad brush strokes follow Derriennic and Lin [4], but with complications. Formally, the solution to the equation $\xi_0 = \eta_0 - B(T)\eta_0$ is $\eta_0 = A(T)\xi_0$, where

$$(10) \qquad A(z) = \frac{1}{1 - B(z)} = \sum_{k=0}^\infty \alpha_k z^k,$$

but there are technicalities in attaching a meaning to $A(T)\xi_0$.

**3. The first step.**

*The size of $R_n$.* The first item of business is to estimate the size of $\|R_n\|$. Here and below, the symbol $\|\cdot\|$ is used more generally to denote the norm in an $L^2$ space, which may vary from one usage to the next.

LEMMA 1.  *Let* $\delta_j = 2^{-j}$. *If* (5) *holds, then*

$$\sum_{j=1}^{\infty} j\sqrt{\ell(2^j)}\sqrt{\delta_j}\|h_{\delta_j}\| < \infty,$$

*where (now)* $\|\cdot\|$ *denotes the norm in* $L^2(\pi)$.

PROOF.  Let $V_n g = g + Qg + \cdots + Q^{n-1}g$, so that $V_n g(w) = E[S_n|W_1 = w]$ and $\|V_n g\| \le 2\|X_0\| + \|E(S_n|\mathcal{F}_0)\|$. Then, rearranging terms in (6),

$$\|h_{\delta_j}\| \le \delta_j \sum_{n=1}^{\infty} \frac{\|V_n g\|}{(1+\delta_j)^n}$$

and

$$\sum_{j=1}^{\infty} j\sqrt{\ell(2^j)}\sqrt{\delta_j}\|h_{\delta_j}\| \le \sum_{n=1}^{\infty}\left[\sum_{j=1}^{\infty} \frac{j\sqrt{\ell(2^j)\delta_j^3}}{(1+\delta_j)^n}\right]\|V_n g\|.$$

Comparing the inner sum to an integral for any fixed integer $n \ge 0$, then

$$\sum_{j=1}^{\infty} \frac{j\sqrt{\ell(2^j)\delta_j^3}}{(1+\delta_j)^n} \le \log_2(e) \int_0^1 \frac{\sqrt{t\ell(2/t)}\log(2/t)}{(1+(1/2)t)^n}\,dt.$$

By a change of variables and the dominated convergence theorem, using Potter's bound (cf. [3], page 25) to supply a dominating function, the integral on the right-hand side of the last inequality is just

$$\frac{1}{\sqrt{n^3}}\int_0^n \sqrt{t\ell\left(\frac{2n}{t}\right)}\log\left(\frac{2n}{t}\right)\left(1+\frac{t}{2n}\right)^{-n}dt \sim \frac{\sqrt{\ell(n)}\log(n)}{\sqrt{n^3}}\int_0^{\infty} \sqrt{t}e^{-(1/2)t}\,dt,$$

from which the lemma follows.  $\square$

PROPOSITION 1.  *If* (5) *holds, then*

$$(11) \qquad \lim_{n\to\infty}\sqrt{\ell(n)}\frac{\|R_n\|}{\sqrt{n}} = 0 \quad and \quad \sum_{n=1}^{\infty}\sqrt{\frac{\ell(n)}{n^3}}\|R_n\| < \infty.$$

PROOF.  Let $H_{\varepsilon}(w_0, w_1) = h_{\varepsilon}(w_1) - Qh_{\varepsilon}(w_0)$, and $M_n(\varepsilon) = H_{\varepsilon}(W_0, W_1) + \cdots + H_{\varepsilon}(W_{n-1}, W_n)$. Then, it is shown in [11] that $S_n = M_n(\varepsilon) + R_n(\varepsilon)$ for each $\varepsilon > 0$ with $R_n(\varepsilon) = \varepsilon S_n(h_{\varepsilon}) + Qh_{\varepsilon}(W_0) - Qh_{\varepsilon}(W_n)$ and $S_n(h_{\varepsilon}) = h_{\varepsilon}(W_1) + \cdots + h_{\varepsilon}(W_n)$. So,

$$R_n = M_n(\varepsilon) - M_n + \varepsilon S_n(h_{\varepsilon}) + Qh_{\varepsilon}(W_0) - Qh_{\varepsilon}(W_n)$$

and

$$\|R_n\| \le \|M_n(\varepsilon) - M_n\| + (n\varepsilon + 2)\|h_{\varepsilon}\| \le \sqrt{n}\|H_{\varepsilon} - H\| + (n\varepsilon + 2)\|h_{\varepsilon}\|.$$

Now let $\varepsilon_n = 2^{-k_n}$, where $2^{k_n-1} \le n < 2^{k_n}$. Then $1/(2n) \le \varepsilon_n = \delta_{k_n} \le 1/n$, and $\|H_{\delta_{j+1}} - H_{\delta_j}\| \le 4\sqrt{\delta_j}\|h_{\delta_j}\|$, by Lemma 2 of [11],

$$\|R_n\| \le \sqrt{n} \sum_{j=k_n}^{\infty} \|H_{\delta_{j+1}} - H_{\delta_j}\| + 3\|h_{\delta_{k_n}}\| \le 10\sqrt{n} \sum_{j=k_n}^{\infty} \sqrt{\delta_j}\|h_{\delta_j}\|.$$

Since $k_n \le j$ implies $n < 2^j$, and so

$$\sum_{k_n \le j} \frac{\sqrt{\ell(n)}}{n} \le \sqrt{\ell(2^j)} \sum_{n < 2^j} \frac{1}{n} \le 2j\sqrt{\ell(2^j)},$$

then we derive

$$\sum_{n=1}^{\infty} \sqrt{\frac{\ell(n)}{n^3}}\|R_n\| \le 10 \sum_{j=1}^{\infty} \left[ \sum_{k_n \le j} \frac{\sqrt{\ell(n)}}{n} \right] \sqrt{\delta_j}\|h_{\delta_j}\|$$

$$\le 20 \sum_{j=1}^{\infty} \sqrt{\ell(2^j)} j \sqrt{\delta_j}\|h_{\delta_j}\|,$$

which is finite by the previous lemma. Thus, the series in (11) converges. That $\sqrt{\ell(n)}\|R_n\|/\sqrt{n} \to 0$ then follows from the subadditivity of $\|R_n\|$; $\|R_{m+n}\| \le \|R_m\| + \|R_n\|$. Since $\|R_n\| \le \|R_k\| + \|R_{n-k}\|$ for all $k = 1, \dots, n-1$, therefore,

$$\sqrt{\frac{\ell(n)}{n}}\|R_n\| \le 6\sqrt{\frac{\ell(n)}{n^3}} \sum_{(1/4)n \le k \le (3/4)n} \|R_k\| \le 6 \sum_{(1/4)n \le k \le (3/4)n} \sqrt{\frac{\ell(k)}{k^3}}\|R_k\|$$

for all sufficiently large $n$, and this approaches 0 as already shown.  □

*The size of $\alpha_n$.* Let

(12)
$$\beta_k = \frac{c}{k} \sum_{n=k}^{\infty} \frac{1}{\sqrt{n^3\ell(n)}}$$

where $c$ is chosen so that $\beta_1 + \beta_2 + \cdots = 1$. Then, $B(z) = \sum_{k=1}^{\infty} \beta_k z^k$ converges for all $|z| \le 1$ in (8) and $\mathcal{R}B(z) < 1$ for all $z \ne 1$, so that $A(z)$ is well defined in (10) for all $|z| \le 1$, except $z = 1$. Observe that $A(z)[1 - B(z)] = 1$ and, therefore,

(13)
$$\alpha_n = \sum_{k=1}^{n} \beta_k \alpha_{n-k}$$

for $n \ge 1$ and $\alpha_0 = 1$. Let

(14)
$$b(t) = B(e^{it}) = \sum_{k=1}^{\infty} \beta_k e^{ikt}$$

for $-\pi < t \le \pi$.

PROPOSITION 2. *b is twice differentiable on $-\pi < t \neq 0 < \pi$, $|1 - b(t)| \sim \kappa_0 \sqrt{|t|}/\sqrt{\ell(1/|t|)}$, and*

$$(15) \qquad |b'(t)| \sim \frac{2c\sqrt{\pi}}{\sqrt{|t|\ell(1/|t|)}}, \qquad |b''(t)| \sim \frac{\kappa_2}{\sqrt{|t|^3\ell(1/|t|)}}$$

*as $t \to 0$, where $\kappa_0 \neq 0$ and $\kappa_2$ are constants (identified) in the proof.*

PROOF. Clearly (14) is absolutely convergent, $b$ is continuous and $b(0) = 1$. By Theorem 2.6 of Zygmund ([19], page 4), the formal expression for the derivative

$$(16) \qquad b'(t) = i \sum_{k=1}^{\infty} \left[ \sum_{n=k}^{\infty} \frac{c}{\sqrt{n^3\ell(n)}} \right] e^{ikt}$$

converges uniformly on $\varepsilon \leq |t| \leq \pi$ for any $\varepsilon > 0$, and therefore, is the derivative of $b$. By Theorem 4.3.2 of [3], page 207,

$$|b'(t)| \sim \frac{2c\sqrt{\pi}}{\sqrt{|t|\ell(1/|t|)}}$$

as $t \to 0$. So, $|1 - b(t)| \sim 4c\sqrt{\pi|t|}/\sqrt{\ell(1/|t|)}$. Reversing the order of summation in (16) (which can be justified by truncating the outer sum at $K$ and letting $K \to \infty$) gives us

$$b'(t) = i \sum_{n=1}^{\infty} \left[ \sum_{k=1}^{n} e^{ikt} \right] \frac{c}{\sqrt{n^3\ell(n)}} = \frac{e^{it}}{1 - e^{it}} \sum_{n=1}^{\infty} (1 - e^{int}) \frac{ic}{\sqrt{n^3\ell(n)}} = f(t)g(t),$$

where $f(t) = e^{it}/(1 - e^{it})$ is continuously differentiable on $-\pi < t \neq 0 < \pi$, and $g$ is continuous. As above,

$$g'(t) = \sum_{n=1}^{\infty} e^{int} \frac{c}{\sqrt{n\ell(n)}}$$

converges uniformly on $\varepsilon \leq |t| \leq \pi$ and

$$|g'(t)| \sim c\sqrt{\pi} \frac{1}{\sqrt{|t|\ell(1/|t|)}}$$

as $t \to 0$. Hence, $b$ is twice continuously differentiable on $-\pi < t \neq 0 < \pi$, and the second relationship in (15) follows from $b''(t) = f'(t)g(t) + f(t)g'(t) = f(t)g'(t) + [ib'(t)/(1 - e^{it})]$ and symmetry. □

In (10), $A(z)$ is defined for all $|z| \leq 1$, except $z = 1$. Let $a(t) = A(e^{it})$ for $-\pi < t \neq 0 < \pi$; then one can derive the following properties.

COROLLARY 3. *a is twice differentiable on $0 < |t| < \pi$, and*

$$|a'(t)| \sim \frac{1}{8c\sqrt{\pi}} \frac{\sqrt{\ell(1/|t|)}}{\sqrt{|t|^3}} \quad and \quad |a''(t)| = O\left(\frac{\sqrt{\ell(1/|t|)}}{\sqrt{|t|^5}}\right)$$

*as $t \to 0$.*

PROOF. This follows directly from (10) and Proposition 2. □

PROPOSITION 3. *Let $\alpha_n$ be the coefficients of $A(z)$; then $0 < \alpha_n \le 1$ for all $n \ge 0$ and*

$$\alpha_n - \alpha_{n+1} = O\left(\frac{\sqrt{\ell(n)}}{\sqrt{n^3}}\right)$$

*as $n \to \infty$.*

PROOF. The first assertion follows easily from (13) and induction. By Proposition 2, $a$ is absolutely integrable, so that $2\pi\alpha_n = \int_{-\pi}^{\pi} e^{-int} a(t)\, dt$, and then

$$\alpha_n - \alpha_{n+1} = \frac{1}{2\pi} \int_{-\pi}^{\pi} e^{-int} a_*(t)\, dt,$$

where $a_*(t) = [1 - e^{-it}]a(t)$. Both $a_*'(s)$ and $sa_*''(s)$ are integrable over $(-\pi, \pi)$. Hence, integration by parts (twice) is justified and yields

$$\alpha_n - \alpha_{n+1} = \frac{1}{2\pi in} \int_{-\pi}^{\pi} e^{-int} a_*'(t)\, dt = \frac{1}{2\pi n^2} \int_{-\pi}^{\pi} [1 - e^{-int}] a_*''(t)\, dt.$$

By Corollary 3, there is a $C$ for which $|a_*''(t)| \le C\sqrt{\ell(1/|t|)/|t|^3}$ for all $0 < |t| \le \pi$. So

$$\begin{aligned}
|\alpha_n - \alpha_{n+1}| &= \frac{1}{2\pi n^3} \left| \int_{-\pi n}^{\pi n} [1 - e^{-it}] a_*''\left(\frac{t}{n}\right) dt \right| \\
&\le \frac{C}{2\pi n^3} \int_{-\pi n}^{\pi n} |1 - e^{-it}| \sqrt{\frac{n^3}{|t|^3} \ell\left(\frac{n}{|t|}\right)}\, dt \\
&\sim \frac{C}{2\pi} \sqrt{\frac{\ell(n)}{n^3}} \int_{-\infty}^{\infty} |1 - e^{-it}| \frac{dt}{\sqrt{|t|^3}},
\end{aligned}$$

using Potter's theorem again and monotonicity of $\ell$. This establishes the proposition. □

*Existence of $\eta_0$.* We need the following fact which is easily deduced from Lemma 1.3 of Krengel ([9], page 4): Let $L_0^2(P)$ be the set of $\eta \in L^2(P)$ with mean 0; if $\theta$ is ergodic, then $[I - T]L_0^2(P)$ is dense in $L_0^2(P)$. Recall the definition of $\xi_0$ in (7) and the expression for $B(T)$ in (9); observe that $\xi_0 \in L_0^2(P)$; and let $A_N(T) = \sum_{n=0}^N \alpha_n T^n$ and $U_n = T + \cdots + T^n$.

PROPOSITION 4. *If* (5) *is satisfied, then* $\eta_0 = \lim_{N \to \infty} A_N(T)\xi_0$ *exists in* $L^2(P)$, *and* $\xi_0 = [I - B(T)]\eta_0$.

PROOF. From (7), we have $U_n \xi_0 = R_n$. Then, summing by parts,

$$A_N(T)\xi_0 = \xi_0 + \alpha_N R_N + \sum_{n=1}^{N-1} (\alpha_n - \alpha_{n+1}) R_n.$$

In view of Propositions 1 and 3 and Karamata's theorem, the sum converges in $L^2(P)$ and $\alpha_N R_N \to 0$.

For the second assertion, let $\eta_N = A_N(T)\xi_0$. Then, rearranging terms and using (13),

$$B(T)\eta_N = \sum_{k=1}^\infty \beta_k \sum_{j=0}^N \alpha_j T^{j+k} \xi_0$$

$$= \sum_{m=1}^N \alpha_m T^m \xi_0 + \sum_{m=N+1}^\infty \left[ \sum_{j=0}^N \alpha_j \beta_{m-j} \right] T^m \xi_0$$

$$= \eta_N - \xi_0 + C_N(T)\xi_0$$

where $C_N(T) := I - [I - B(T)]A_N(T)$. So, it suffices to show that $\|C_N(T)\xi_0\| \to 0$. For this, first observe that, replacing $T$ by $z$ in the definition of $C_N(T)$, $1 - C_N(z) = [1 - B(z)]A_N(z)$. Then $C_N(1) = 1$ and the coefficients of $C_N(z)$ are all positive, so that $\|C_N(T)\|_{\mathrm{op}} \le 1$, where $\|\cdot\|_{\mathrm{op}}$ stands for operator norm. So, it suffices to show that $\|C_N(T)\xi\| \to 0$ for all $\xi \in [I - T]L_0^2(P)$, a dense subset of $L_0^2(P)$. This is easy: for if $\xi = \psi - T\psi$, then

$$C_N(T)\xi = \sum_{j=0}^N \alpha_j \left[ \beta_{N+1-j} T^{N+1} \psi + \sum_{m=N+1}^\infty (\beta_{m+1-j} - \beta_{m-j}) T_m \psi \right]$$

and

$$\|C_N(T)\xi\| \le 2\|\psi\| \sum_{j=0}^N \alpha_j \beta_{N+1-j} \to 0$$

as $N \to \infty$ by (13) and Proposition 3. $\square$

**4. The second step.** Some preparation is necessary for the second step. First, for any $\eta \in L^2(P)$, $\eta^* := \sup_{n \geq 1} U_n |\eta|/n \in L^2(P)$ by the dominated ergodic theorem (see, e.g., Krengel [9], page 52). We will also use the following fact:

$$(17) \qquad\qquad E(\sqrt{(\eta^2)^*}) \leq 2\|\eta\|,$$

whose proof is essentially an application of the maximal ergodic theorem ([14], Corollary 2.2) to $(\eta^2)^*$.

The proof of Theorem 1 will be completed by proving:

THEOREM 2. *If* $\xi \in [I - B(T)]L^2(P)$, *then*

$$\lim_{n \to \infty} \frac{U_n \xi}{\sqrt{n \ell^*(n)}} = 0 \qquad w.p.\ 1.$$

PROOF. By assumption, there is an $\eta \in L^2(P)$ for which $\xi = \eta - B(T)\eta = \sum_{k=1}^{\infty} \beta_k[\eta - T^k \eta]$, and there is no loss of generality in supposing that $\eta \in L_0^2(P)$. Observe that $|T^k \eta|^p = T^k(|\eta|^p)$ for any integer $k \geq 0$ and real $p > 0$, and write

$$U_n \xi = I_n \eta + II_n \eta,$$

where

$$I_n \eta = \sum_{k=1}^{n} \beta_k U_n[\eta - T^k \eta]$$

and

$$II_n \eta = \sum_{k=n+1}^{\infty} \beta_k U_n[\eta - T^k \eta].$$

If $k > n$, then $|U_n(\eta - T^k \eta)| \leq |U_n \eta| + |U_n T^k \eta| \leq [\eta^* + T^k \eta^*]n$. So,

$$|II_n \eta| \leq n \sum_{k=n+1}^{\infty} \beta_k[\eta^* + T^k \eta^*].$$

Here

$$\sum_{k=n+1}^{\infty} \beta_k T^k \eta^* \leq \sum_{k=n+1}^{\infty} \Delta \beta_k U_k \eta^* \leq \sum_{k=n+1}^{\infty} k \Delta \beta_k \eta^{**},$$

where $\Delta \beta_k = \beta_k - \beta_{k+1}$ and $\eta^{**} = \sup_{k \geq 1} U_k \eta^*/k$. Observing that

$$\sum_{k=n+1}^{\infty} (\beta_k + k \Delta \beta_k) = n \beta_{n+1} + 2 \sum_{k=n+1}^{\infty} \beta_k,$$

thus,

$$|II_n\eta| \leq n(\eta^* \vee \eta^{**})\left[\sum_{k=n+1}^{\infty}\beta_k + \sum_{k=n+1}^{\infty}k\Delta\beta_k\right] = (\eta^* \vee \eta^{**}) \times O\left(\sqrt{\frac{n}{\ell(n)}}\right)$$

and

$$(18) \qquad \lim_{n\to\infty}\frac{II_n\eta}{\sqrt{n\ell^*(n)}} = 0 \qquad \text{w.p. 1.}$$

Similarly, for $k \leq n$, $U_n\eta - U_nT^k\eta = U_k\eta - U_kT^n\eta$; then

$$I_n\eta = \sum_{k=1}^{n}\beta_kU_k\eta - \sum_{k=1}^{n}\beta_kU_kT^n\eta.$$

Letting $\gamma_j = \sum_{k=j}^{\infty}\beta_k$ and recalling (12), we have

$$\sum_{j=1}^{n}\gamma_j^2 \sim (4c)^2\left(\sum_{j=1}^{n}\frac{1}{j\ell(j)}\right) = (4c)^2\ell^*(n)$$

and

$$|I_n\eta| \leq \sum_{k=1}^{n}\beta_k\sum_{j=1}^{k}[T^j|\eta| + T^{j+n}|\eta|] \leq \sum_{j=1}^{n}\gamma_j[T^j|\eta| + T^{j+n}|\eta|]$$

$$\leq \sqrt{\sum_{j=1}^{n}\gamma_j^2} \times \sqrt{2 \times \sum_{j=1}^{2n}T^j\eta^2}.$$

Using (17), there exists a constant $C > 0$, such that

$$E\left(\sup_n\frac{|I_n\eta|}{\sqrt{n\ell^*(n)}}\right) \leq C\|\eta\|,$$

where $C$ does not depend on $\eta$. Hence, to show

$$(19) \qquad \lim_{n\to\infty}\frac{I_n\eta}{\sqrt{n\ell^*(n)}} = 0 \qquad \text{w.p. 1}$$

for each $\eta \in L_0^2(P)$, one only needs to consider $\eta \in (I - T)L_0^2(P)$, a dense subset in $L_0^2(P)$, and this is easy. If $\eta = \phi - T\phi$ for some $\phi \in L_0^2(P)$, then $U_kT^n\eta = T^{n+1}\phi - T^{k+n+1}\phi$ for $1 \leq k \leq n$, so that

$$|I_n\eta| \leq \left|T\sum_{k=1}^{n}\beta_k(\phi - T^k\phi)\right| + \left|T^{n+1}\sum_{k=1}^{n}\beta_k(\phi - T^k\phi)\right| \leq T\tilde{\phi} + T^{n+1}\tilde{\phi},$$

where

$$\tilde{\phi} = \sum_{k=1}^{\infty}\beta_k|\phi - T^k\phi| \in L^2(P).$$

Since $\tilde{\phi} \in L^2(P)$, $\lim_{n\to\infty} T^{n+1}\tilde{\phi}/\sqrt{n} = 0$ w.p. 1 by an easy application of the Borel–Cantelli Lemma and therefore, $\lim_{n\to\infty} I_n \eta/\sqrt{n\ell^*(n)} = 0$ w.p. 1. The theorem now follows by combining (18) and (19).  □

**5. Invariance principles.** Let $C[0,1]$ be the space of all real-valued continuous functions on $[0,1]$, endowed with the metric

$$\rho(x,y) = \sup_{0 \le t \le 1} |x(t) - y(t)|,$$

where $x, y \in C[0,1]$. For any $\nu \ge 0$, let $K_\nu$ denote the set of absolutely continuous functions $x \in C[0,1]$ such that $x(0) = 0$ and

$$\int_0^1 [x'(t)]^2 \, dt \le \nu^2.$$

Set $S_0 = M_0 = 0$ and define sequences of random functions $\{\theta_n(\cdot)\}$ and $\{\zeta_n(\cdot)\}$ respectively by

$$\theta_n(t) = \frac{S_k + (nt - k)X_{k+1}}{\sqrt{2n \log_2(n)}},$$

$$\zeta_n(t) = \frac{M_k + (nt - k)(M_{k+1} - M_k)}{\sqrt{2n \log_2(n)}},$$

for $k \le nt \le k+1, k = 0, 1, \ldots, n-1$. Then $\theta_n, \zeta_n \in C[0,1]$.

COROLLARY 4. *If the hypothesis in Corollary* 1 *holds, then w.p.* 1, $\{\theta_n\}_{n \ge 3}$ *are relatively compact in* $C[0,1]$, *and the set of limit points is* $K_\sigma$.

PROOF. Under the hypothesis, (3) and (4) hold; then

$$\rho(\theta_n, \zeta_n) \le \max_{k \le n} \frac{|R_k|}{\sqrt{2n \log_2(n)}} \to 0 \qquad \text{w.p. 1,}$$

which implies that $\theta_n$ and $\zeta_n$ have the same limit points; and the limit points of $\zeta_n$ are known to be $K_\sigma$ w.p. 1 (see, e.g., Heyde and Scott [7], Corollary 2). □

Let

$$\mathbb{B}_n(t) = \frac{1}{\sqrt{n}} S_{\lfloor nt \rfloor}$$

for $0 \le t < 1$, $\mathbb{B}_n(1) = \mathbb{B}_n(1-)$, where $\lfloor \cdot \rfloor$ denotes the integer part. Then $\mathbb{B}_n \in D[0,1]$, the space of càdlàg functions as described in Chapter 3 of Billingsley [2]. Let $F_n$ denote a regular conditional distribution for $\mathbb{B}_n$ given $\mathcal{F}_0$, so that $F_n(\omega; B) = P[\mathbb{B}_n \in B | \mathcal{F}_0](\omega)$ for Borel sets $B \subseteq D[0,1]$; and let $\Phi_\sigma$ denote the distribution of $\sigma\mathbb{B}$, where $\mathbb{B}$ is a standard Brownian motion. Let $\Delta$ denote the Prokhorov metric on $D[0,1]$ (cf. [2], page 238).

COROLLARY 5.   *If the hypothesis in Corollary* 2 *holds, then*

$$\lim_{n \to \infty} \Delta[F_n(\omega; \cdot), \Phi_\sigma] = 0 \qquad a.e. \ \omega. \tag{20}$$

PROOF.   For $S_n = M_n + R_n$, let $M_n^*(t) = M_{\lfloor nt \rfloor}/\sqrt{n}$, $0 \leq t < 1$ and $M_n^*(1) = M_n^*(1-)$. Let $G_n$ denote a regular conditional distribution for the random element $M_n^*$ given $\mathcal{F}_0$. Then $G_n(\omega; \cdot)$ converges to $\Phi_\sigma$ for a.e. $\omega$ $(P)$, by verifying Theorem 2.5 of Durrett and Resnick [5] in view of the mean ergodic theorem. Under the hypothesis of Corollary 2, $\max_{1 \leq k \leq n} |R_k|/\sqrt{n} \to 0$ w.p. 1, and therefore,

$$\rho(M_n^*, \mathbb{B}_n) = \sup_{0 \leq t \leq 1} |M_n^*(t) - \mathbb{B}_n(t)| \to 0 \qquad \text{w.p. 1.}$$

Equation (20) follows.   □

**6. Examples.**   In this section, we illustrate our conditions by considering linear processes, additive functionals of a Bernoulli shift and $\rho$-mixing processes.

*Linear processes.*   Let $\dots, \varepsilon_{-1}, \varepsilon_0, \varepsilon_1, \dots$ be an ergodic stationary martingale difference sequence with common mean 0 and variance 1. Define a linear process

$$X_k = \sum_{j=0}^{\infty} a_j \varepsilon_{k-j},$$

where $a_0, a_1, \dots$ is a square summable sequence, and observe that $X_k$ is of the form $g(W_k)$ with $W_k = (\dots, \varepsilon_{k-1}, \varepsilon_k)$.

PROPOSITION 5.   *Suppose* $a_n = O[1/(nL(n))]$, *where* $L(\cdot)$ *is a positive, nondecreasing, slowly varying function. If*

$$\sum_{n=2}^{\infty} \frac{\log^\alpha(n)}{nL(n)} < \infty \tag{21}$$

*with* $\alpha = 3/2$, *then* (5) *holds with* $\ell(n) = 1 \vee \log(n)$ *and, thus the conclusions to Corollaries* 1 *and* 4. *Furthermore, if* (21) *holds with some* $\alpha > 3/2$, *then also the conclusions to Corollaries* 2 *and* 5 *hold.*

PROOF.   Letting $s_{j,n} = a_{j+1} + \dots + a_{j+n}$, straightforward calculations yield that

$$\|E[S_n | \mathcal{F}_0]\|^2 = \sum_{j=0}^{\infty} s_{j,n}^2.$$

If $j \geq 3$, then

$$|s_{j,n}| \leq \frac{C}{L(j)} \int_j^{j+n} \frac{1}{x} \, dx \leq \frac{C}{L(j)} \log\left(1 + \frac{n}{j}\right)$$

for some constant $C > 0$, and therefore,

$$\sum_{j=3}^{\infty} s_{j,n}^2 \leq C^2 \int_2^{\infty} \frac{1}{L^2(x)} \log^2\left(1 + \frac{n}{x}\right) dx$$

$$= nC^2 \int_0^{n/2} \frac{1}{L^2(n/t)} \frac{\log^2(1+t)}{t^2} \, dt = O\left[\frac{n}{L^2(n)}\right],$$

where the last step follows from the dominated convergence theorem, using Potter's bound to supply the dominating function, or by Fatou's lemma. It is then easily verified that $\|E(S_n|\mathcal{F}_0)\| = O[\sqrt{n}/L(n)]$, and the proposition is an immediate consequence. $\square$

REMARK 1. If $L(n) \sim \log^{\beta}(n)$, then (21) requires $\beta > 5/2$. This is similar to, but not strictly comparable with, the results of Yokoyama [18], who required finite moments of order $p > 2$ and $\beta \geq 1 + (2/p)$.

*Additive functionals of the Bernoulli shift.* Now consider a Bernoulli process, say

$$W_k = \sum_{j=1}^{\infty} \frac{1}{2^j} \varepsilon_{k-j+1},$$

where $\ldots, \varepsilon_{-1}, \varepsilon_0, \varepsilon_1, \ldots$ are i.i.d. random variables that take the values 0 and 1 with probability $1/2$ each. Then $\mathcal{W} = [0, 1]$, $\pi$ is the uniform distribution, and

$$Qf(w) = \frac{1}{2}\left[f\left(\frac{w}{2}\right) + f\left(\frac{1+w}{2}\right)\right]$$

for $f \in L^1$. Next, consider a stationary process of the form $X_k = g(W_k)$, where $g$ is square integrable with respect to $\pi$ and has mean 0. In this case, it is possible to relate (5) to a weak regularity condition on $g$.

PROPOSITION 6. *If*

$$(22) \qquad \int_0^1 \int_0^1 \frac{[g(x) - g(y)]^2}{|x - y|} \log^{5/2+\delta}\left[\log\left(\frac{1}{|x-y|}\right)\right] dx \, dy < \infty$$

*for some $\delta > 0$, then the conclusions to Corollaries 2 and 5 hold, and so also those of Corollaries 1 and 4.*

PROOF (Sketched). The proof involves showing that (22) implies (5), for which $\ell(n)$ can be chosen such that $\ell^*(n)$ remains bounded. The details are similar to the proof of Proposition 3 in [11], and will be omitted. $\square$

*ρ-mixing processes.* Our condition (5) can be checked when a mixing rate is available for a *ρ*-mixing process; see [12], pages 4–5 for a definition.

COROLLARY 6. *Let* $\rho(n)$ *be the ρ-mixing coefficients of a centered, square integrable, stationary process* $(X_k)_{k \in \mathbb{Z}}$. *If* $\rho(n) = O(\log^\gamma n)$ *for some* $\gamma > 5/2$, *as* $n \to \infty$, *then* (1) *holds.*

PROOF (Outline). Let $S_n = X_1 + \cdots + X_n$ and $h(x) = (1 \vee \log x)^{3/2}$. By an argument similar to that in [12], page 15, one can easily show that, for some constant $C > 0$,

$$\sum_{r=0}^\infty \frac{h(2^r)\|E(S_{2^r}|\mathcal{F}_0)\|}{2^{r/2}} \leq C \sum_{j=0}^\infty h(2^j)\rho(2^j) < \infty.$$

Since $\|E(S_n|\mathcal{F}_0)\|$ is subadditive, it is then straightforward to argue as in Lemma 2.7 of [13], that

$$\sum_{n=1}^\infty \frac{h(n)\|E(S_n|\mathcal{F}_0)\|}{n^{3/2}} < \infty.$$

Therefore, (1) holds by Corollary 1. □

REMARK 2. Shao [16] showed that LIL holds when $\rho(n) = O(\log^\gamma n)$ for some $\gamma > 1$, but through a completely different approach.

**Acknowledgments.** Thanks to a referee and an Associate Editor for their useful suggestions, including the application to *ρ*-mixing sequences and the current proof of Theorem 2 using (17). An earlier one had used Banach's principle. Thanks also to Dalibor Volny and Karl Petersen for helpful conversations.

DEPARTMENT OF STATISTICS                     DEPARTMENT OF STATISTICS
UNIVERSITY OF MICHIGAN                        UNIVERSITY OF MICHIGAN
439 WEST HALL                                462 WEST HALL
ANN ARBOR, MICHIGAN                           ANN ARBOR, MICHIGAN
USA                                          USA
E-MAIL: ouzhao@umich.edu                     E-MAIL: michaelw@umich.edu